  \documentclass[11pt]{amsart}
\usepackage[dvips]{graphicx}
\usepackage{textcomp}
\usepackage{amsbsy}
\usepackage{latexsym}
\usepackage[mathscr]{eucal}
\usepackage{amsfonts,amsmath,amsthm}
\usepackage{amssymb}
\usepackage[usenames]{color}
\usepackage{fullpage}

\begin{document}
\bibliographystyle{plain}
 \title{Rescaled Pure  Greedy Algorithm for Hilbert and Banach Spaces}
 \author{Guergana Petrova}
 \thanks{%
 This research was supported by the Office of
 Naval Research ContractÊ Ê ONR N00014-11-1-0712, and
 byÊ the NSF  Grant DMS 1222715.}%
\hbadness=10000
\vbadness=10000
\newtheorem{lemma}{Lemma}[section]
\newtheorem{prop}[lemma]{Proposition}
\newtheorem{cor}[lemma]{Corollary}
\newtheorem{theorem}[lemma]{Theorem}
\newtheorem{remark}[lemma]{Remark}
\newtheorem{example}[lemma]{Example}
\newtheorem{definition}[lemma]{Definition}
\newtheorem{proper}[lemma]{Properties}
\newtheorem{assumption}[lemma]{Assumption}
\newenvironment{disarray}{\everymath{\displaystyle\everymath{}}\array}{\endarray}

\def\RR{\rm \hbox{I\kern-.2em\hbox{R}}}
\def\NN{\rm \hbox{I\kern-.2em\hbox{N}}}
\def\ZZ{\rm {{\rm Z}\kern-.28em{\rm Z}}}
\def\CC{\rm \hbox{C\kern -.5em {\raise .32ex \hbox{$\scriptscriptstyle
|$}}\kern
-.22em{\raise .6ex \hbox{$\scriptscriptstyle |$}}\kern .4em}}
\def\vp{\varphi}
\def\<{\langle}
\def\>{\rangle}
\def\t{\tilde}
\def\i{\infty}
\def\e{\varepsilon}
\def\sm{\setminus}
\def\nl{\newline}
\def\o{\overline}
\def\wt{\widetilde}
\def\wh{\widehat}
\def\cT{{\cal T}}
\def\cA{{\cal A}}
\def\cI{{\cal I}}
\def\cV{{\cal V}}
\def\cB{{\cal B}}
\def\cF{{\cal F}}

\def\cR{{\cal R}}
\def\cD{{\cal D}}
\def\cP{{\cal P}}
\def\cJ{{\cal J}}
\def\cM{{\cal M}}
\def\cO{{\cal O}}
\def\Chi{\raise .3ex
\hbox{\large $\chi$}} \def\vp{\varphi}
\def\lsima{\hbox{\kern -.6em\raisebox{-1ex}{$~\stackrel{\textstyle<}{\sim}~$}}\kern -.4em}
\def\lsim{\hbox{\kern -.2em\raisebox{-1ex}{$~\stackrel{\textstyle<}{\sim}~$}}\kern -.2em}
\def\[{\Bigl [}
\def\]{\Bigr ]}
\def\({\Bigl (}
\def\){\Bigr )}
\def\[{\Bigl [}
\def\]{\Bigr ]}
\def\({\Bigl (}
\def\){\Bigr )}
\def\L{\pounds}
\def\pr{{\rm Prob}}
\newcommand{\cs}[1]{{\color{magenta}{#1}}}
\def\ds{\displaystyle}
\def\ev#1{\vec{#1}}     
\newcommand{\lt}{\ell^{2}(\nabla)}
\def\Supp#1{{\rm supp\,}{#1}}
\def\R{\mathbb{R}}
\def\E{\mathbb{E}}
\def\nl{\newline}
\def\T{{\relax\ifmmode I\!\!\hspace{-1pt}T\else$I\!\!\hspace{-1pt}T$\fi}}
\def\N{\mathbb{N}}
\def\Z{\mathbb{Z}}
\def\N{\mathbb{N}}
\def\Zd{\Z^d}
\def\Q{\mathbb{Q}}
\def\C{\mathbb{C}}
\def\Rd{\R^d}
\def\gsim{\mathrel{\raisebox{-4pt}{$\stackrel{\textstyle>}{\sim}$}}}
\def\sime{\raisebox{0ex}{$~\stackrel{\textstyle\sim}{=}~$}}
\def\lsim{\raisebox{-1ex}{$~\stackrel{\textstyle<}{\sim}~$}}
\def\div{\mbox{ div }}
\def\M{M}  \def\NN{N}                  
\def\L{{\ell}}               
\def\Le{{\ell^1}}            
\def\Lz{{\ell^2}}
\def\Let{{\tilde\ell^1}}     
\def\Lzt{{\tilde\ell^2}}
\def\Ltw{\ell^\tau^w(\nabla)}
\def\t#1{\tilde{#1}}
\def\la{\lambda}
\def\La{\Lambda}
\def\ga{\gamma}
\def\BV{{\rm BV}}
\def\Ga{\eta}
\def\al{\alpha}
\def\cZ{{\cal Z}}
\def\cA{{\cal A}}
\def\cU{{\cal U}}
\def\argmin{\mathop{\rm argmin}}
\def\argmax{\mathop{\rm argmax}}
\def\prob{\mathop{\rm prob}}

\def\cO{{\cal O}}
\def\cA{{\cal A}}
\def\cC{{\cal C}}
\def\cS{{\cal F}}
\def\bu{{\bf u}}
\def\bz{{\bf z}}
\def\bZ{{\bf Z}}
\def\bI{{\bf I}}
\def\cE{{\cal E}}
\def\cD{{\cal D}}
\def\cG{{\cal G}}
\def\cI{{\cal I}}
\def\cJ{{\cal J}}
\def\cM{{\cal M}}
\def\cN{{\cal N}}
\def\cT{{\cal T}}
\def\cU{{\cal U}}
\def\cV{{\cal V}}
\def\cW{{\cal W}}
\def\cL{{\cal L}}
\def\cB{{\cal B}}
\def\cG{{\cal G}}
\def\cK{{\cal K}}
\def\cX{{\cal X}}
\def\cS{{\cal S}}
\def\cP{{\cal P}}
\def\cQ{{\cal Q}}
\def\cR{{\cal R}}
\def\cU{{\cal U}}
\def\bL{{\bf L}}
\def\bl{{\bf l}}
\def\bK{{\bf K}}
\def\bC{{\bf C}}
\def\X{X\in\{L,R\}}
\def\ph{{\varphi}}
\def\D{{\Delta}}
\def\H{{\cal H}}
\def\bM{{\bf M}}
\def\bx{{\bf x}}
\def\bj{{\bf j}}
\def\bG{{\bf G}}
\def\bP{{\bf P}}
\def\bW{{\bf W}}
\def\bT{{\bf T}}
\def\bV{{\bf V}}
\def\bv{{\bf v}}
\def\bt{{\bf t}}
\def\bz{{\bf z}}
\def\bw{{\bf w}}
\def \span{{\rm span}}
\def \meas {{\rm meas}}
\def\rhom{{\rho^m}}
\def\diff{\hbox{\tiny $\Delta$}}
\def\EE{{\rm Exp}}
\def\lan{\langle}
\def\ran{\rangle}
\def\argmin{\mathop{\rm argmin}}
\def\codim{\mathop{\rm codim}}
\def\rank{\mathop{\rm rank}}

\def\argmax{\mathop{\rm argmax}}
\def\dJ{\nabla}
\newcommand{\ba}{{\bf a}}
\newcommand{\bb}{{\bf b}}
\newcommand{\bc}{{\bf c}}
\newcommand{\bd}{{\bf d}}
\newcommand{\bs}{{\bf s}}
\newcommand{\bff}{{\bf f}}
\newcommand{\bp}{{\bf p}}
\newcommand{\bg}{{\bf g}}
\newcommand{\by}{{\bf y}}
\newcommand{\br}{{\bf r}}
\newcommand{\be}{\begin{equation}}
\newcommand{\ee}{\end{equation}}
\newcommand{\bea}{$$ \begin{array}{lll}}
\newcommand{\eea}{\end{array} $$}
\def \Vol{\mathop{\rm  Vol}}
\def \mes{\mathop{\rm mes}}
\def \Prob{\mathop{\rm  Prob}}
\def \exp{\mathop{\rm    exp}}
\def \sign{\mathop{\rm   sign}}
\def \sp{\mathop{\rm   span}}
\def \rad{\mathop{\rm   rad}}
\def \vphi{{\varphi}}
\def \csp{\overline \mathop{\rm   span}}
%
%
\newcommand{\beqn}{\begin{equation}}
\newcommand{\eeqn}{\end{equation}}
\def\beginproof{\noindent{\bf Proof:}~ }
\def\endproof{\hfill\rule{1.5mm}{1.5mm}\\[2mm]}

\newenvironment{Proof}{\noindent{\bf Proof:}\quad}{\endproof}

\renewcommand{\theequation}{\thesection.\arabic{equation}}
\renewcommand{\thefigure}{\thesection.\arabic{figure}}

\makeatletter
\@addtoreset{equation}{section}
\makeatother

\newcommand\abs[1]{\left|#1\right|}
\newcommand\clos{\mathop{\rm clos}\nolimits}
\newcommand\trunc{\mathop{\rm trunc}\nolimits}
\renewcommand\d{d}
\newcommand\dd{d}
\newcommand\diag{\mathop{\rm diag}}
\newcommand\dist{\mathop{\rm dist}}
\newcommand\diam{\mathop{\rm diam}}
\newcommand\cond{\mathop{\rm cond}\nolimits}
\newcommand\eref[1]{{\rm (\ref{#1})}}
\newcommand{\iref}[1]{{\rm (\ref{#1})}}
\newcommand\Hnorm[1]{\norm{#1}_{H^s([0,1])}}
\def\int{\intop\limits}
\renewcommand\labelenumi{(\roman{enumi})}
\newcommand\lnorm[1]{\norm{#1}_{\ell^2(\Z)}}
\newcommand\Lnorm[1]{\norm{#1}_{L_2([0,1])}}
\newcommand\LR{{L_2(\R)}}
\newcommand\LRnorm[1]{\norm{#1}_\LR}
\newcommand\Matrix[2]{\hphantom{#1}_#2#1}
\newcommand\norm[1]{\left\|#1\right\|}
\newcommand\ogauss[1]{\left\lceil#1\right\rceil}
\newcommand{\QED}{\hfill
\raisebox{-2pt}{\rule{5.6pt}{8pt}\rule{4pt}{0pt}}%
  \smallskip\par}
\newcommand\Rscalar[1]{\scalar{#1}_\R}
\newcommand\scalar[1]{\left(#1\right)}
\newcommand\Scalar[1]{\scalar{#1}_{[0,1]}}
\newcommand\Span{\mathop{\rm span}}
\newcommand\supp{\mathop{\rm supp}}
\newcommand\ugauss[1]{\left\lfloor#1\right\rfloor}
\newcommand\with{\, : \,}
\newcommand\Null{{\bf 0}}
\newcommand\bA{{\bf A}}
\newcommand\bB{{\bf B}}
\newcommand\bR{{\bf R}}
\newcommand\bD{{\bf D}}
\newcommand\bE{{\bf E}}
\newcommand\bF{{\bf F}}
\newcommand\bH{{\bf H}}
\newcommand\bU{{\bf U}}
\newcommand\cH{{\cal H}}
\newcommand\sinc{{\rm sinc}}
\def\enorm#1{| \! | \! | #1 | \! | \! |}

\newcommand{\dm}{\frac{d-1}{d}}

\let\bm\bf
\newcommand{\bbeta}{{\mbox{\boldmath$\beta$}}}
\newcommand{\bal}{{\mbox{\boldmath$\alpha$}}}
\newcommand{\bbi}{{\bm i}}

\def\nnew{\color{Red}}
\def\mnew{\color{Blue}}

\newcommand{\dI}{\Delta}
\newcommand\aconv{\mathop{\rm absconv}}

\maketitle
\date{}

 \DeclareGraphicsRule{.tif}{png}{.png}{`convert #1 `dirname
 #1`/`basename #1 .tif`.png}
 
  \begin{abstract}
 We 
 show that a very simple modification of the Pure  Greedy Algorithm for approximating functions  by sparse sums from a dictionary in a Hilbert or more generally a Banach space has 
  optimal convergence rates on the class of convex combinations of dictionary elements.\\
  
  \end{abstract}

 \section{Introduction}
 \label{Intr}
 Greedy algorithms have been used quite extensively as a tool for generating   approximations from  redundant families of functions, such as  frames or more general dictionaries  
 ${\mathcal D}$.   Given a Banach space $X$,    a dictionary  is any set  ${\mathcal D}$ of norm one elements from $  X$ whose span is dense in $X$.
  The most natural greedy algorithm in a Hilbert space is the Pure Greedy Algorithm ({\bf PGA}), which is  also known as Matching Pursuit, see \cite{DT}  for the description of this and other algorithms.    The fact that  the {\bf PGA}  lacks optimal convergence properties has led to a variety of modified greedy algorithms such as the Relaxed Greedy Algorithm
  ({\bf RGA}), the Orthogonal Greedy Algorithm, and 
 their weak versions. There are also analogues of these, developed for approximating functions in Banach spaces, see \cite{Tbook}. 
 
 The  central issues in the study of these algorithms is their ease of implementation and    their approximation power, measured in terms of convergence rates. If 
 $f_m$ is the output of a greedy algorithm after $m$ iterations,  then $f_m$ is a linear combination of at most $m$ dictionary elements.  Such linear combinations are said to be sparse of order $m$.
 The quality of the approximation is measured by  the decay of the  error $\|f-f_m\|$ as $m\rightarrow \infty$, where $\|\cdot\|$ is the norm in the Hilbert or Banach space, respectively.  Of course, the decay rate of this error is governed by properties of the target function $f$.   The typical properties imposed on $f$ are that it is sparse, or more generally, that it is in some way compressible.  Here, compressible    means that it can be written as a (generally speaking, infinite) linear combination of dictionary elements with some restrictions on the coefficients.   The most
  frequently applied assumption on $f$ is that it is in the unit ball of the class ${\mathcal A}_1({\mathcal D})$, that is  the set of all functions which are a convex combination of dictionary elements
  (provided we consider symmetric dictionaries).  It is known that the elements in this class can be approximated by $m$ sparse vectors to accuracy ${\mathcal O}(m^{-1/2})$, see Theorem \ref{cr},  and so this rate of approximation serves as a benchmark for the performance of greedy algorithms.
  
 It has been shown in \cite{DT} in the case of Hilbert space that  whenever $f\in {\mathcal A}_1({\mathcal D})$,  the output $f_m$ of the {\bf PGA} satisfies
 \be
 \label{pga}
 \|f-f_m\|={\mathcal O}(m^{-1/6}),\quad m\to\infty.
 \ee
 Later results gave slight improvements of the  above estimate. For example, in \cite{KT},  the rate ${\mathcal O}(m^{-1/6})$ was improved to  ${\mathcal O}(m^{-11/62})$. Based on the method from the latter paper, 
Sil'nichenko  \cite{S} then showed   a rate   of ${\mathcal O}(m^{-\frac{s}{2(s+2)}})$, where $s$ solves a certain equation, and that $\frac{s}{2(s+2)}>11/62$. Similar estimates  for the weak versions of the {\bf PGA} can be found in \cite{Tbook}.
Estimates for the error from below have also been provided, see \cite{LT,L}.

The fact that the {\bf PGA} does not attain the optimal rate for approximating the elements in ${\mathcal A}_1({\mathcal D})$ has led to various modifications of this algorithm.   Two of these modifications, the  Relaxed and the Orthogonal Greedy Algorithm  were shown to achieve the optimal rate ${\mathcal O}(m^{-1/2})$, see \cite{DT}. 
 
 The purpose of the present paper is to show that a very simple modification of the {\bf PGA}, namely just rescaling $f_m$ at each iteration, already leads to the improved convergence rate ${\mathcal O}(m^{-1/2})$ for functions in ${\mathcal A}_1({\mathcal D})$. 
  The rescaling we suggest is simply the orthogonal projection of 
$f$ onto $f_m$.     
We call this modified algorithm a Rescaled Pure Greedy Algorithm  ({\bf RPGA})  and  prove optimal convergence rates for its weak version in Hilbert and Banach spaces.
In a subsequent paper, see \cite {GP}, we show that this strategy can also be applied successfully for developing an algorithm for convex optimization.

The paper is organized as follows.  In \S{\ref{cond}, we spell out our notation and recall some simple known facts related to greedy algorithms. In \S\ref{H}, we  present the {\bf RPGA}   for a Hilbert space and prove the above convergence rate. The  remaining parts of this paper consider a
modification of this algorithm for Banach spaces and weak versions of this algorithm.

 \section{ Notation and Preliminaries}
 \label{cond}
  
   We denote by $H$ a Hilbert space and by $X$ a Banach space with $\|\cdot\|$  being the norm in 
 these spaces, respectively. A set of functions
 ${\mathcal D}\subset H(\hbox{or}\,\,X)$ is called a dictionary 
 if $\|\varphi\|=1$ for every $\varphi\in {\mathcal D}$
 and the closure of $span {(\mathcal D)}$ is $H(\hbox{or}\,\,X)$. An example of
 a dictionary is any Shauder basis for $H(\hbox{or}\,\,X)$.  However, the main
 idea behind dictionaries is to cover redundant families such as frames.  A common example of dictionaries is the
   union of several Shauder bases. 
  
 The set $\Sigma_m({\mathcal D})$  consists of all $m$-sparse elements with respect to the dictionary ${\mathcal D}$, namely
 $$
 \Sigma_m:= \Sigma_m({\mathcal D})=\{g: \,\,g=\sum_{\varphi\in\Lambda} c_{\varphi}\varphi, \,\,\Lambda\in {\mathcal D},\, \,|\Lambda|\leq m\}.
 $$
Here, we use the notation $|\Lambda|$ to denote the cardinality of the index set $\Lambda$.  For a general element $f$ from  $X$, we define   the error of approximation  
   $$
 \sigma_m(f):=\sigma_m(f,{\mathcal D}):=\inf_{g\in \Sigma_m}\|f-g\|
 $$
of $f$  by elements from $\Sigma_m$. 
The rate of decay of $\sigma_m(f)$ as $m\to \infty$ says how well $f$ can be approximated by sparse elements.     
 
  For a general dictionary ${\mathcal D}\subset H(\hbox{or}\,\,X)$,  we define the class of functions 
  $$
 {\mathcal A}^{o}_1({\mathcal D}, M):=\{f=\sum_{k \in \Lambda} c_k(f)\varphi_k:  \, \varphi_k\in{\mathcal D}, \,|\Lambda|<\infty, \, 
 \sum_{k \in \Lambda} |c_k(f)|\leq M\},
 $$
 and by ${\mathcal A}_1({\mathcal D}, M)$ its closure in $H(\hbox{or}\,\,X)$.   Then,  ${\mathcal A}_1({\mathcal D})$ is defined to be the union of the classes 
 ${\mathcal A}_1({\mathcal D}, M)$ over all $M>0$. For $f\in {\mathcal A}_1({\mathcal D})$,  we define the ``semi-norm'' of $f$ as
 $$
 |f|_{{\mathcal A}_1(\mathcal D)}:=\inf \{M:\,\,f\in {\mathcal A}_1({\mathcal D}, M)\}.
 $$

 A fundamental result for approximating ${\mathcal A}_1({\mathcal D})$ is the following, see \cite{DT}.
 \begin{theorem}
 \label{cr}
For a general dictionary  ${\mathcal D}\subset H$
and $f\in {\mathcal A}_1({\mathcal D})\subset H$, we have
 $$
 \sigma_m(f,{\mathcal D})\leq c|f|_{{\mathcal A}_1({\mathcal D})}m^{-1/2}, \quad  m=1,2,\ldots.
 $$
\end{theorem}

 When analyzing the convergence of greedy algorithms, we will use  the following lemma,
 proved in \cite{NP}. 
 \begin{lemma}
 \label{lmseq}
 Let $\ell>0$, $r>0$, $B>0$, and 
 $\{a_m\}_{m=1}^{\infty}$ andÊ $\{r_m\}_{m=2}^{\infty}$
 be sequences of non-negative numbers satisfying the
 inequalities
 $$ a_1\leq B, \quad a_{m+1} \leq a_m(1-
 \frac{r_{m+1}}{r}a_m^\ell), \quad m=1,2,\dots.$$
 Then, we have 
 \begin{equation}
 \label{tuti12}
 a_m \leq
 \max\{1,\ell^{-1/\ell}\}r^{1/\ell}(rB^{-\ell}+\Sigma_{k=2}^m
 r_k)^{-1/\ell}, \quad m=2,3, \ldots.
 \end{equation}
 \end{lemma}
 We note that several similar versions of this lemma have been proved and used in analysis of greedy algorithms, see \cite{Tbook}.
 \section{The Hilbert space case}
 \label{H}
In  order to show the simplicity of our results, we  begin with the standard case of the {\bf RPGA} in a Hilbert space.   Later, we treat the case of Banach spaces and weak algorithms, but the reader familiar with this topic will see that the results in these more general settings follow by standard modifications of the results from this section.    We denote the inner product in the Hilbert space $H$ by $\langle\cdot,\cdot\rangle$, and so
 the norm of $f\in H$ is  $\|f\|=\langle f,f\rangle^{1/2}$.
 
  The {\bf  RPGA($ \mathcal D$)} is defined by the following simple  steps.
 \bigskip
 
\noindent {\bf RPGA($ \mathcal D$):}
   \begin{itemize}
 \item {\bf Step $0$}: 
 Define $f_0:=0$.  
 \item {\bf Step $m$}: 
 \item Assuming $f_{m-1}$ has been
 computed and $f_{m-1}\neq f$. Choose a direction  $\varphi_m\in {\mathcal D}$ such that
 $$
 |\langle f-f_{m-1},\varphi_m\rangle|=\sup_{\varphi\in{\mathcal D}} |\langle f-f_{m-1},\varphi\rangle|.
 $$
With 
$$\lambda_m:=\langle f-f_{m-1},\varphi_m\rangle,
 \quad 
 \hat f_m:=f_{m-1}+\lambda_m\varphi_{m},
 \quad 
 s_m:=\frac{\lan f, \hat f_m\ran}{\|\hat f_m\|^2},
  $$
   define the next approximant to be
 $$
 f_m=s_m\hat f_m.
 $$
 \item If $f=f_{m}$, stop the algorithm and define $f_k=f_{m}=f$, for $k>m$.
 \item If $f\neq f_m$,  proceed to Step $m+1$.
 \end{itemize}
 \bigskip
\noindent    
 Note that if the output at each Step $m$ were $\hat f_m$ and not $f_m=s_m\hat f_m$,   this would be  the {\bf PGA}. 
However, the new algorithm uses not $\hat f_m$, but the best approximation to $f$ from the one dimensional space $span\{\hat f_m\}$, that is
 $s_m\hat f_m$.   Adding this step,
which is just appropriate scaling of the output of the {\bf PGA}, allows us to prove optimal convergence rate of $m^{-1/2}$ for the 
proposed algorithm.

Next, we show that the {\bf RPGA} and the Relaxed Greedy Algorithm ({\bf RGA}) provide different sequences of approximants
$\{f_m\}$ and $\{f_m^r\}$, respectively, and thus {\bf RPGA} is different from the known so far greedy algorithms. For both algorithms
$$
f_0=f_0^r=0, \quad f_1=f_1^r=\langle f,\varphi_1\rangle \varphi_1, 
$$
where $\varphi_1\in{\mathcal D}$ is such that  $ |\langle f,\varphi_1\rangle|=\sup_{\varphi\in{\mathcal D}} |\langle f,\varphi\rangle|$. For both
\bf RPGA} and {\bf RGA}, the next element $\varphi_2\in{\mathcal D}$ is chosen as
$|\langle f-f_{1},\varphi_2\rangle|=\sup_{\varphi\in{\mathcal D}} |\langle f-f_{1},\varphi\rangle|$. One can easily compute that the 
next approximant, generated by the {\bf RPGA} is
$$
f_2=s_2f_1+
s_2 \langle f-f_1,\varphi_2\rangle\varphi_2, \quad s_2=\frac{\langle f,\varphi_1\rangle^2+\langle f,\varphi_2\rangle^2
-\langle f,\varphi_1\rangle\langle f,\varphi_2\rangle\langle \varphi_1,\varphi_2\rangle}{
\langle f,\varphi_1\rangle^2+\langle f,\varphi_2\rangle^2
-\langle f,\varphi_1\rangle^2\langle \varphi_1,\varphi_2\rangle^2},
$$
while the classical  {\bf RGA} would give
$$
f_2^r=\frac{1}{2}f_1+\frac{1}{2}\varphi_2.
$$
There are some modifications of the {\bf RGA},  see \cite{B}, where the approximant at Step $m$  is determined not as
$$
f_m^r=(1-\frac{1}{m})f_{m-1}^r+\frac{1}{m}\varphi_m, \quad  \hbox{where}\,\,\,|\langle f-f_{m-1}^r,\varphi_m\rangle|=\sup_{\varphi\in{\mathcal D}} |\langle f-f_{m-1}^r,\varphi\rangle|,
$$
but as
\begin{equation}
\label{pl}
f^r_m=(1-a_m)f^r_{m-1}+a_m\varphi_m, 
\end{equation}
where $a_m$ and $\varphi_m$ are the solutions of the minimization problem
$$
\min_{a\in [0,1], \varphi\in {\mathcal D}}\|f-\left((1-a)f^r_{m-1}+a\varphi\right)\|.
$$
While the sequence, generated by the {\bf RPGA}  is a linear combination of $f_{m-1}$ and $\varphi_m$, that is
$$
f_m=s_mf_{m-1}+\lambda_ms_m\varphi_m,
$$
it is different from the convex combinations \eref{pl}, from other variations of the {\bf RGA}, as described in \cite{Tbook}, and from 
the best approximation to $f$ from $span\{f_{m-1}^r,\varphi_m\}$. For example, the best approximation to 
$f$ from $span\{f_{1}^r,\varphi_2\}$ is
$$
f_2^r=\frac{\langle f,\varphi_1\rangle^2-\langle f,\varphi_1\rangle\langle f,\varphi_2\rangle\langle \varphi_1,\varphi_2\rangle}
{\langle f,\varphi_1\rangle^2(1-\langle \varphi_1,\varphi_2\rangle^2)}f_1+\frac{\langle f,\varphi_2\rangle
-\langle f,\varphi_1\rangle\langle \varphi_1,\varphi_2\rangle}{1-\langle \varphi_1,\varphi_2\rangle^2}\varphi_2,
$$
and again $f_2\neq f_2^r$.
In summary, we can view the new algorithm either as a rescaled version of the {\bf PGA} or a new modification  of the {\bf RGA}.

We continue with the following theorem.

\begin{theorem}
\label{HSconv}
If $f\in{\mathcal A}_1({\mathcal D})\subset H$, then the output $(f_m)_{m\ge 0}$  of the  {\bf RPGA}($ \mathcal D$)  satisfies
\be
\label{optdga}
e_m:=\|f-f_m\|\le |f|_{{\mathcal A}_1({\mathcal D})}m^{-1/2}, \quad m=1,2\dots.
\ee
\end{theorem}

 \noindent{\bf Proof:}  Since $f_m$ is the orthogonal projection of $f$ onto the one dimensional space spanned by $\hat f_m$, we have 
 \be
 \label{HSorthogonality}
 \langle f-f_m,f_m\rangle=0, \quad m\ge 0.
 \ee
  Next, note that
   the definition of $\hat f_m$ and the choice of $\lambda_m$  give
\begin{eqnarray}
\nonumber
\|f-\hat f_m\|^2&=& \langle f-f_{m-1}-\lambda_m\varphi_m, f-f_{m-1}-\lambda_m\varphi_m\rangle
 \nonumber \\
 &=&
 \|f-\hat f_{m-1}\|^2-2\lambda_m \langle f-f_{m-1},\varphi_m\rangle+\lambda_m^2\|\varphi_m\|^2
 \nonumber \\
 &=&
 \|f-f_{m-1}\|^2-\langle f-f_{m-1},\varphi_m\rangle^2,
 \label{popa}
\end{eqnarray}
where we have used that $\|\varphi_m\|=1$.   Now, assume  $f\neq f_{m-1}$. Since $f_m$ 
 is the orthogonal projection of $f$ onto $span\{\hat f_{m}\}$, we have
  $$
 e_m^2=\|f-f_m\|^2=\|f-s_m\hat f_m\|^2\leq \|f-\hat f_m\|^2.
 $$
 We combine the latter inequality and \eref{popa}  to derive that 
 \begin{equation}
 \label{poi}
 e_m^2\leq e_{m-1}^2-\langle f-f_{m-1},\varphi_m\rangle^2, \quad m=1,2,\ldots.
 \end{equation}
We proceed with an estimate from below for $\langle f-f_{m-1},\varphi_m\rangle$.
 Note that 
 \begin{equation}
 \label{rim}
 e_{m-1}^2=\|f-f_{m-1}\|^2=\langle f-f_{m-1},f-f_{m-1}\rangle=\langle f-f_{m-1},f\rangle,
 \end{equation}
 where we have used \eref{HSorthogonality}.
 
 It is enough to prove \eref{optdga} for functions $f$ that are finite sums $f=\sum_jc_j\varphi_j$ with $\sum_j|c_j|\leq M$, since these functions are dense 
 in ${\mathcal A}_1({\mathcal D},M)$.
 Let us fix $\varepsilon >0$ and choose a representation
 for $f=\sum_{\varphi\in {\mathcal D}}c_\varphi \varphi$,
 such that 
 $$
\sum_{\varphi\in {\mathcal D}}|c_\varphi|<M+\varepsilon.
 $$
 It follows from \eref{rim} that 
 \begin{eqnarray}
 \nonumber
 e_{m-1}^2&=&\sum_{\varphi\in {\mathcal D}}c_\varphi\lan
f-f_{m-1},\varphi\ran\\
 \nonumber
 &\leq&  |\lan f-f_{m-1},\varphi_{m}\ran|\sum_{\varphi\in {\mathcal D}}|c_\varphi|
 \nonumber \\
 &<& |\lan f-f_{m-1},\varphi_{m}\ran|(M+\varepsilon),
 \nonumber
 \end{eqnarray}
 where we have used the choice of $\varphi_{m}$.
We let  $\varepsilon \rightarrow 0$ and  obtain the inequality
 \begin{equation}
 \label{pop}
 M^{-1}e_{m-1}^2\leq |\lan f-f_{m-1},\varphi_m\ran|.
 \end{equation}
 We combine \eref{poi} and \eref{pop} to obtain
 $$
 e_m^2\leq e_{m-1}^2- M^{-2}e_{m-1}^{4}=e_{m-1}^2(1-M^{-2}e_{m-1}^2),\quad m\ge 2.
 $$
 Note that 
 $$
 \|f\|^2=\lan f,f\ran=\sum_{\varphi\in {\mathcal D}}c_\varphi \lan f,\varphi\ran\leq
 |\lan f,\varphi_1\ran|\sum_{\varphi\in {\mathcal D}}|c_\varphi|<\|f\|(M+\varepsilon), 
 $$
 and therefore $\|f\|\leq M$.
 Since $e_1^2\le e_0^2=\|f\|^2\le M^2$, we can apply Lemma \ref{lmseq} with $a_m=e_m^2$, $B=M^2$, $r_m:=1$, 
$r=M^2$, and $\ell=1$.  Then, \eref{tuti12} gives
 \be
 \nonumber
 e_m^2\le M^2m^{-1}, \quad m\ge 2,
 \ee
 and the theorem follows.
 \hfill $\Box$

 In the sections that follow, we introduce variants of the {\bf RPGA} and prove convergence results
 similar to Theorem \ref{HSconv}.

 \section{The Weak Rescaled Pure  Greedy Algorithm for Hilbert spaces}
 \label{greedy}
 In this section, we describe the Weak Rescaled Pure 
 Greedy Algorithm ({\bf WRPGA}). It is determined by a weakness
 sequence $\{t_k\}_{k=1}^\infty$,
 where all $t_k\in (0,1]$,  and the 
 dictionary $\mathcal D$. We denote it by {\bf WRPGA}($\{t_k\}, \mathcal D$).
 
 \bigskip
 
 \noindent
 {\bf WRPGA}($\{t_k\},  \mathcal D$):
 \begin{itemize}
 \itemÊ {\bf Step $0$}: 
 Define $f_0=0$.Ê  
 \item {\bf Step $m$}:Ê 
 \item Assuming $f_{m-1}$ has been
 computed and $f_{m-1}\neq f$. Choose a direction  $\varphi_m\in {\mathcal D}$ such that
 $$
 |\langle f-f_{m-1},\varphi_m\rangle|\geq t_m\sup_{\varphi\in{\mathcal D}} |\langle f-f_{m-1},\varphi\rangle|.
 $$
 With
 $$
\lambda_m=\langle f-f_{m-1},\varphi_m\rangle,
 \quad 
 \hat f_m:=f_{m-1}+\lambda_m\varphi_{m},
 \quad 
 s_m=\frac{\lan f, \hat f_m\ran}{\|\hat f_m\|^2},
  $$
  define the next approximant to be
 $$
 f_m=s_m\hat f_m.
 $$
 \item If $f=f_{m-1}$, stop the algorithm and define $f_k=f_{m-1}=f$ for $k\geq m$.
  \item If $f\neq f_m$,  proceed to Step $m+1$.
 \end{itemize}
 \bigskip
\noindent    
In the case when all elements $t_k$ of the weakness sequence are $t_k=1$, this algorithm is  the {\bf RPGA}($\mathcal D$).   The following theorem 
holds.
\begin{theorem}
\label{WHSconv}
If $f\in{\mathcal A}_1({\mathcal D})\subset H$, then the output $(f_m)_{m\ge 0}$  of  the {\bf WRPGA}($ \{t_k\},\mathcal D$)  satisfies
\begin{eqnarray}
\label{woptdga}
e_m:=\|f-f_m\|\le |f|_{{\mathcal A}_1({\mathcal D})}\left(\sum_{k=1}^mt_k^2\right)^{-1/2}, \quad m\ge 1.
\end{eqnarray}
\end{theorem}
 \noindent{\bf Proof:}  The proof is similar to the one of Theorem \ref{HSconv}, where
 we show that for  the error $e_m^2=\|f-f_m\|^2$, we have the inequality,
 \begin{equation}
 \label{wpoi}
 e_m^2\leq e_{m-1}^2-\langle f-f_{m-1},\varphi_m\rangle^2, \quad m=1,2,\ldots.
 \end{equation}
The estimate from below for $\langle f-f_{m-1},\varphi_m\rangle$
 is derived similarly as
 \begin{equation}
 \label{wpop}
 M^{-1}t_me_{m-1}^2\leq |\lan f-f_{m-1},\varphi_m\ran|,
 \end{equation}
 where we have used the definition of $\varphi_m$.
 Next, it follows from  \eref{wpoi} and \eref{wpop} that
 $$
 e_m^2\leq e_{m-1}^2- M^{-2}t_m^2e_{m-1}^{4}=e_{m-1}^2(1-M^{-2}t_m^2e_{m-1}^2),\quad m\ge 1,
 $$
 Note that 
 $$
 \|f\|^2=\lan f,f\ran=\sum_{\varphi\in {\mathcal D}}c_\varphi \lan f,\varphi\ran\leq
t_1^{-1} |\lan f,\varphi_1\ran|\sum_{\varphi\in {\mathcal D}}|c_\varphi|<t_1^{-1}\|f\|(M+\varepsilon), 
 $$
 and therefore $\|f\|\leq Mt_1^{-1}$.
 Since $e_1^2\le e_0^2=\|f\|^2\le M^2t_1^{-2}$, we can apply Lemma \ref{lmseq} with $a_m=e_m^2$, $B=M^2t_1^{-2}$, $r_m:=t_m^2$, 
$r=M^2$, and $\ell=1$ to obtain
 \be
 \nonumber
 e_m^2\le M^2\left(t_1^2+\sum_{k=2}^mt_k^2\right)^{-1}, \quad m\ge 2,
 \ee
 and the theorem follows.
 \hfill $\Box$

 \section{The Banach space case}
 \label{B}
 \noindent
  In this section, we will state the  {\bf RPGA}($\mathcal D$) algorithm for Banach spaces $X$  with norm $\|\cdot\|$ and dictionary $\mathcal D$, and prove 
  convergence results for certain Banach spaces. Let us first start with the introduction of the modulus of smoothness $\rho$ of a Banach
  space $X$ , which is defined as
  $$
  \rho(u):=\sup_{f,g\in X, \|f\|=\|g\|=1}\left \{\frac{1}{2}(\|f+ug\|+\|f-ug\|)-1 \right \}, \quad u>0.
  $$
 In this paper, we shall consider only Banach spaces $X$ whose modulus of smoothness satisfies the inequality
  $$
  \rho(u)\leq \gamma u^q, \quad 1<q\leq 2,\quad \gamma\hbox{\,\,-constant}.
  $$
This is a natural assumption, since the modulus of smoothness of $X=L_p$, $1<p<\infty$, for example, is known to satisfy such inequality.  Recall that, see \cite{DG}, 
for $X=L_p$, 
$$ \rho(u)\leq
\left
\{\begin{array}{cc}
\displaystyle{
\frac{1}{p}u^p}, & \mbox{     if  } 1\leq p\leq 2, \\\\
\displaystyle{\frac{p-1}{2}u^2}, & \mbox{     if  }  2\leq p<\infty.
\end{array}
\right.
$$

 Next, for every element $f\in X$, $f\neq 0$, we consider its  norming functional $F_f\in X^*$ with the properties
  $ \|F_f\|=1$, $ F_f(f)=\|f\|$.
  Note that if $X=H$ is a Hilbert space, the norming functional for $f\in H$ is
  $$
  F_f(\cdot)=\frac{<f,\cdot>}{\|f\|}.
  $$
  There is a relationship between the norming functional  $F_g$ for any $g\in X$, $g\neq 0$,  and the modulus of smoothness of $X$, given by the
  following lemma.
  \begin{lemma}
  \label{qq}
  Let $X$ be a Banach space with modulus of smoothness $\rho$, where $\rho(u)\leq \gamma u^q$,  $1<q\leq 2$. Let  
  $g\in X$, $g\neq 0$ with norming functional $F_g$. Then, for every $h\in X$ , we have
  \begin{equation}
  \|g+uh\|\leq \|g\|+uF_g(h)+2\gamma u^q\|g\|^{1-q}\|h\|^q, \quad u>0.
  \end{equation}
  \end{lemma}
  \noindent{\bf Proof:}
  The proof follows from Lemma 6.1 in \cite{Tbook} and the property of the modulus of smoothness.
   \hfill $\Box$

   We next present the {\bf RPGA}($\mathcal D$) for the Banach space $X$ with dictionary $\mathcal D$.
   \bigskip
 
 \noindent
 {\bf RPGA}($\mathcal D$):
 \begin{itemize}
 \itemÊ {\bf Step $0$}: 
 Define $f_0=0$.Ê  
 \item {\bf Step $m$}:Ê 
 \item Assuming $f_{m-1}$ has been
 computed and $f\neq f_{m-1}$. Choose a direction  $\varphi_m\in {\mathcal D}$ such that
 $$
|F_{ f-f_{m-1}}(\varphi_m)|=\sup_{\varphi\in{\mathcal D}}|F_{ f-f_{m-1}}(\varphi)|.
 $$
With 
 $$
\lambda_m=sign\{F_{ f-f_{m-1}}(\varphi_m)\}\|f-f_{m-1}\|(2\gamma q)^{\frac{1}{1-q}}|F_{ f-f_{m-1}}(\varphi_m)|^{\frac{1}{q-1}},
\quad
 \hat f_m:=f_{m-1}+\lambda_m\varphi_{m},
$$
choose $s_m$ such that 
$$
\|f- s_m\hat f_m\|=\min_{s\in \R}\|f- s\hat f_m\|,
  $$
  and define the next approximant to be
 $$
 f_m=s_m\hat f_m.
 $$
\item  If $f=f_{m-1}$, stop the algorithm and define $f_k=f_{m-1}=f$ for $k\geq m$.
\item If $f\neq f_m$,  proceed to Step $m+1$.
 \end{itemize}
\noindent

The  following lemma holds.
 \begin{lemma}
  \label{qq1}
  Let $X$ be a Banach space with modulus of smoothness  $\rho$, $\rho(u)\leq \gamma u^q$,  $1<q\leq 2$.
  Let $f_{m-1}$ be the output of the {\bf RPGA}($ \mathcal D$) at Step $m-1$. Then, if $f\neq f_{m-1}$, we
  have 
  $$
  F_{f-f_{m-1}}(f_{m-1})=0.
  $$
  \end{lemma}
  \noindent{\bf Proof:}    
 Let us denote by $L:=span\{\hat f_{m-1}\}\subset X$. Clearly,  $f_{m-1}\in L$, and moreover, $f_{m-1}$ is the best 
 approximation to $f$ from $L$. We apply Lemma 6.9 from \cite{Tbook} to the linear space $L$ and the vector $f_{m-1}$, and derive the
 lemma.
\hfil$\Box$

The next theorem provides the convergence rate for the new algorithm in Banach spaces.

\begin{theorem}
\label{XSconv}
Let $X$ be a Banach space with modulus of smoothness $\rho(u)\leq \gamma u^q$, 
$1<q\leq 2$. If $f\in{\mathcal A}_1({\mathcal D})\subset X$, then the output $(f_m)_{m\ge 0}$  of  the {\bf RPGA}($ \mathcal D$)  satisfies
\be
\label{woptdga}
e_m:=\|f-f_m\|\le c |f|_{{\mathcal A}_1({\mathcal D})}m^{1/q-1}, \quad m\ge 2,
\ee
where $c=c(\gamma,q)$.
\end{theorem}

 \noindent{\bf Proof:}    
 Clearly, we haveÊ $e_0=\|f-f_0\|=\|f\|$.
At Step $m$, $m=1,2,\ldots$  of the algorithm, either $f=f_{m-1}$, in which case $f_k=f_{m-1}$, $k\geq m$, and therefore $e_m=0$, or we have
\begin{eqnarray}
\nonumber
e_m&=&\|f-f_m\|=\|f-s_m\hat f_m\|\leq \|f-\hat f_m\|
=\|(f-f_{m-1})-\lambda_m\varphi_m\|.
\nonumber 
\end{eqnarray}
We now apply Lemma \ref{qq} to the latter inequality with $g=f-f_{m-1}\neq 0$, $u=|\lambda_m|>0$,  $h=-sign\{\lambda_m\}\varphi_m$, and derive
\begin{eqnarray}
\label {az}
e_m&\leq& \|f-f_{m-1}\|-\lambda_mF_{f-f_{m-1}}(\varphi_m)+2\gamma |\lambda_m|^q\|f-f_{m-1}\|^{1-q}\|\varphi_m\|^q
\nonumber  \\
&=& e_{m-1}-\lambda_mF_{f-f_{m-1}}(\varphi_m)+2\gamma |\lambda_m|^qe_{m-1}^{1-q}
\nonumber  \\
&=& e_{m-1}-\frac{q-1}{q}\left (2\gamma q\right)^{\frac{1}{1-q}}e_{m-1}|F_{f-f_{m-1}}(\varphi_m)|^{\frac{q}{q-1}},
\end{eqnarray}
where we have used that $\|\varphi_m\|=1$ and the choice of $\lambda_m$. Now, we need an estimate from below for $|F_{f-f_{m-1}}(\varphi_m)|$.
Using Lemma \ref{qq1}, we obtain that 
\begin{eqnarray}
\label{pol}
e_{m-1}=\|f-f_{m-1}\|=F_{f-f_{m-1}}(f-f_{m-1})=F_{f-f_{m-1}}(f).
\end{eqnarray}
As in the Hilbert space case, it  is enough to consider  functions $f$ that are finite sums $f=\sum_jc_j\varphi_j$ with $\sum_j|c_j|\leq M$, since these functions are dense 
 in ${\mathcal A}_1({\mathcal D},M)$.
 Let us fix $\varepsilon >0$ and choose a representation
 for $f=\sum_{\varphi\in {\mathcal D}}c_\varphi \varphi$,
 such that 
 $$
\sum_{\varphi\in {\mathcal D}}|c_\varphi|<M+\varepsilon.
$$
 It follows that 
\begin{eqnarray}
\nonumber
 F_{f-f_{m-1}}(f)&=&\sum_{\varphi\in {\mathcal D}}c_\varphi F_{f-f_{m-1}}(\varphi)\leq \sum_{\varphi\in {\mathcal D}}|c_\varphi|| F_{f-f_{m-1}}(\varphi)|
 \nonumber \\
 &\leq& | F_{f-f_{m-1}}(\varphi_m)|\sum_{\varphi\in {\mathcal D}}|c_\varphi|
< | F_{f-f_{m-1}}(\varphi_m)|(M+\varepsilon).
 \nonumber
\end{eqnarray}
We take $\epsilon \rightarrow 0$ and derive
$$
F_{f-f_{m-1}}(f)\leq | F_{f-f_{m-1}}(\varphi_m)|M
$$
The latter estimate and \eref{pol} provide the estimate from below
$$
M^{-1}e_{m-1}\leq | F_{f-f_{m-1}}(\varphi_m)|,
$$
which together with \eref{az} result in
$$
e_m\leq e_{m-1}\left (1-\frac{q-1}{q}\left (2\gamma q\right)^{\frac{1}{1-q}}M^{-\frac{q}{q-1}}e_{m-1}^{\frac{q}{q-1}}\right).
$$
Note that $e_1\leq e_0=\|f\|\leq M$, since
$$
\|f\|=F_f(f)=\sum_\varphi c_\varphi F_f(\varphi)\leq |F_f(\varphi_1)|\sum_\varphi |c_\varphi |<M+\varepsilon,
$$
for every $\varepsilon>0$. We now use Lemma  \ref{lmseq} with $a_m=e_m$, $B=M$, $r_m:=\frac{q-1}{q}\left (2\gamma q\right)^{\frac{1}{1-q}}$, 
$r=M^{\frac{q}{q-1}}$, and $\ell=\frac{q}{q-1}$ to obtain
 \be
 \nonumber
 e_m\le M\left(1+\frac{q-1}{q}\left (2\gamma q\right)^{\frac{1}{1-q}}(m-1)\right)^{1/q-1}, \quad m\ge 2,
 \ee
 and the theorem follows.

\hfill $\Box$

   \section{The Weak Rescaled Pure  Greedy Algorithm for Banach spaces}
 \label{greedy}
 In this section, weÊ describe the Weak Rescaled Pure
 Greedy Algorithm for Banach spaces. It is determined by a weakness
 sequence $\{t_k\}_{k=1}^\infty$,
 where all $t_k\in (0,1]$,  and the 
 dictionary $\mathcal D$. As in the Hilbert case, we  denote it by {\bf WRPGA}($\{t_k\}, \mathcal D$).
 
   \bigskip
 
 \noindent
  {\bf WRPGA}($\{t_k\}, \mathcal D$):
 \begin{itemize}
 \itemÊ {\bf Step $0$}: 
 Define $f_0=0$.Ê  
 \item {\bf Step $m$}:Ê 
 \item Assuming $f_{m-1}$ has been
 computed and $f\neq f_{m-1}$. Choose a direction  $\varphi_m\in {\mathcal D}$ such that
 $$
|F_{ f-f_{m-1}}(\varphi_m)|\geq t_m\sup_{\varphi\in{\mathcal D}}|F_{ f-f_{m-1}}(\varphi)|.
 $$
With 
 $$
\lambda_m=sign\{F_{ f-f_{m-1}}(\varphi_m)\}\|f-f_{m-1}\|(2\gamma q)^{\frac{1}{1-q}}|F_{ f-f_{m-1}}(\varphi_m)|^{\frac{1}{q-1}},
\quad
 \hat f_m:=f_{m-1}+\lambda_m\varphi_{m},
$$
choose $s_m$ such that 
$$
\|f- s_m\hat f_m\|=\min_{s\in \R}\|f- s\hat f_m\|,
  $$
  and define the next approximant to be
 $$
 f_m=s_m\hat f_m.
 $$
\item  If $f=f_{m-1}$, stop the algorithm and define $f_k=f_{m-1}=f$ for $k\geq m$.
\item If $f\neq f_m$,  proceed to Step $m+1$.
 \end{itemize}
\noindent
 Next, we present the convergence rates for the {\bf WRPGA}($\{t_k\}, \mathcal D$) in Banach Spaces.
\begin{theorem}
\label{WXSconv}
Let $X$ be a Banach space with modulus of smoothness $\rho(u)\leq \gamma u^q$, 
$1<q\leq 2$. If $f\in{\mathcal A}_1({\mathcal D})\subset X$, then the output $(f_m)_{m\ge 0}$  of  the {\bf WRPGA}($\{t_k\}, \mathcal D$)  satisfies
\begin{eqnarray}
\label{woptdga}
e_m:=\|f-f_m\|\le c |f|_{{\mathcal A}_1({\mathcal D})}\left(\sum_{k=1}^mt_{k}^{\frac{q}{q-1}}\right)^{1/q-1}, \quad m\geq 1,
\end{eqnarray}
where $c=c(\gamma,q)$.
\end{theorem}

 \noindent{\bf Proof:}    
As in the proof of Theorem \ref{XSconv}, we show that 
\begin{eqnarray}
\label {waz}
e_m\leq e_{m-1}-\frac{q-1}{q}\left (2\gamma q\right)^{\frac{1}{1-q}}e_{m-1}|F_{f-f_{m-1}}(\varphi_m)|^{\frac{q}{q-1}}.
\end{eqnarray}
Next, similarly to  Theorem \ref{XSconv}, we prove an estimate from below for $|F_{f-f_{m-1}}(\varphi_m)|$, which is
$$
M^{-1}t_me_{m-1}\leq | F_{f-f_{m-1}}(\varphi_m)|,
$$
which together with \eref{waz} result in
$$
e_m\leq e_{m-1}\left (1-\frac{q-1}{q}\left (2\gamma q\right)^{\frac{1}{1-q}}t_{m}^{\frac{q}{q-1}}M^{-\frac{q}{q-1}}e_{m-1}^{\frac{q}{q-1}}\right).
$$
Again, since  $e_1\leq e_0=\|f\|\leq Mt_1^{-1}$, we can use Lemma  \ref{lmseq} with $a_m=e_m$, $B=Mt_1^{-1}$, $r_m:=\frac{q-1}{q}\left (2\gamma q\right)^{\frac{1}{1-q}}t_{m}^{\frac{q}{q-1}}$, 
$r=M^{\frac{q}{q-1}}$, and $\ell=\frac{q}{q-1}$.  Then, \eref{tuti12} gives
 \be
 \nonumber
 e_m\le M\left(t_1^\frac{q}{q-1}+\frac{q-1}{q}\left (2\gamma q\right)^{\frac{1}{1-q}}\sum_{k=2}^mt_{k}^{\frac{q}{q-1}}\right)^{1/q-1}, \quad m\ge 2,
 \ee
 and the theorem follows.
\hfill $\Box$  
 
 \vskip .1in
 Ê 
 Ê \noindent
 Guergana Petrova\\
 Department of Mathematics, Texas A\&M University,
 College Station, TX 77843, USA\\
 Ê gpetrova@math.tamu.edu
 
 \end{document}